\documentclass[12pt]{article}
\usepackage{amsmath}
\usepackage{latexsym}
\usepackage{amssymb}
\newtheorem{thm}{Theorem}[section]

\newtheorem{Defn}[thm]{Definition}
\newtheorem{Remark}[thm]{Remark}
\newtheorem{Note}[thm]{Note}

\newtheorem{Example}[thm]{Example}
\newtheorem{Examples}[thm]{Examples}
\newtheorem{Problems}[thm]{Problems}

\newtheorem{Problem}[thm]{Problem}
\newtheorem{Convention}[thm]{Convention}
\newtheorem{Number}[thm]{\!\!}

                  {\nopagebreak\hspace*{\fill}$\Box$\medskip\medskip\par}   
\newcommand{\Punkt}{\nopagebreak\hspace*{\fill}$\Box$}

\newcommand{\at}{\symbol{'100}}

\newcommand{\N}{{\mathbb N}}
\newcommand{\R}{{\mathbb R}}

\newcommand{\C}{{\mathbb C}}

\newcommand{\cW}{{\cal W}}

\newcommand{\sub}{\subseteq}

\newcommand{\cA}{{\cal A}}

\begin{document}
$\;$\\[-24mm]
\begin{center}
{\Large\bf Topological algebras of rapidly decreasing\vspace{2mm} matrices and generalizations}\\[6mm]
{\bf Helge Gl\"{o}ckner and Bastian Langkamp\footnote{Research
supported by DFG grant GL 357/5--1}}\vspace{4mm}
\end{center}
\begin{abstract}\noindent
It is well-known fact in K-theory that the rapidly decreasing matrices of
countable
size form a locally m-convex associative topological algebra whose set
of quasi-invertible elements
is open, and such that the quasi-inversion map is continuous.
We generalize these conclusions to further
algebras of weighted matrices with entries in a Banach algebra.\vspace{3mm}
\end{abstract}
{\footnotesize {\em Classification}:
Primary 46H20;
Secondary 46A45, 22E65\\
{\em Key words}: Rapidly decreasing matrix, weighted matrix algebra, continuous inverse\linebreak
algebra, $Q$-algebra}\\[5mm]
If $(\cA,\|.\|)$ is a Banach algebra over $\R$ or $\C$ and $\cW$ a non-empty set of monotonically
increasing functions $f\colon \N\to\;]0,\infty[$,
we define $M(\cA,\cW)$
as the set of all $T=(t_{ij})_{i,j\in \N}\in \cA^{\N\times\N}$ such that
\[
\|T\|_f:=\sup_{i,j\in \N}\, f(i\vee j)\|t_{ij}\|\,<\,\infty
\]
for all $f\in \cW$, where $i\vee j$ denotes the maximum of $i$ and $j$.
It is clear that $M(\cA,\cW)$ is a vector space;
we give it the locally convex Hausdorff vector topology defined
by the set of norms $\{\|.\|_f\colon f\in \cW\}$. We show:\\[6mm]
{\bf Theorem.}
\emph{Assume there exists $g\in \cW$ such that $C_g:=\sum_{n=1}^\infty \frac{1}{g(n)}<\infty$.
If $R=(r_{ij})_{i,j\in\N}, S=(s_{ij})_{i,j\in\N}\in M(\cA,\cW)$
and $i,j\in\N$, then the series}
\[
t_{ij}:=\sum_{k=1}^\infty r_{ik}s_{kj}
\]
\emph{converges absolutely in $\cA$.
Moreover, $RS:=(t_{ij})_{i,j\in\N}\in M(\cA,\cW)$,
and the multiplication defined in this way makes $M(\cA,\cW)$
a locally m-convex, associative topological algebra
which is complete as a topological vector space,
has an open set of quasi-invertible elements,
and whose quasi-inversion map is continuous.}\vspace{2mm}\pagebreak

\noindent
Recall that an element $x$ in an associative (not necessarily unital) algebra $A$
is called \emph{quasi-invertible} if there exists $y\in A$ such that $xy=yx$
and $x+y-xy=0$. The element $q(x):=y$ is then unique and is called
the \emph{quasi-inverse} of $x$.
Locally convex topological algebras $A$ with an open set $Q(A)$ of quasi-invertible elements
and continuous quasi-inversion map $q\colon Q(A)\to A$
are called \emph{continuous quasi-inverse algebras}
(and \emph{continuous inverse algebras} if they have, moreover, a unit element).
See \cite{Wal} for information on such algebras
as well as \cite{Glo} and \cite{Nee},
where such algebras are inspected due to their usefulness in
infinite-dimensional Lie theory.
Also recall that a topological algebra $A$ is called \emph{locally m-convex}
if its vector topology can be defined using a set of seminorms $p\colon A\to[0,\infty[$
which are sub-multiplicative, i.e., $p(xy)\leq p(x)p(y)$ for all $x,y\in A$.
If $A$ is, moreover, complete as a topological vector space,
this means that $A$ is a projective limit of Banach algebras~\cite{Mic}.\\[3mm]
If we take $\cW:=\{f_m\colon m\in \N_0\}$ with $f_m(n):=n^m$,
and $\cA:=\C$, then $M(\C,\cW)$ is the so-called
algebra of rapidly decreasing
matrices, which plays an important role
in the K-theory of Fr\'{e}chet algebras~\cite{Phi}.
It is known that this algebra (and its counterpart for general $A$)
has an open group of quasi-invertible elements~\cite[4.6]{Phi}
and is a locally m-convex Fr\'{e}chet algebra \cite[2.4\,(1)]{Phi}.
Our discussion recovers these facts,
but applies to larger classes of weighted matrix algebras.
As we realized, only a simple condition
(the existence of $g$ with $C_g<\infty$) needs to be
imposed on the set of weights.\\[5mm]
{\bf Proof of the theorem.}
Step 1. Let $R=(r_{ij})_{i,j\in \N}$
and $S=(s_{ij})_{i,j\in \N}$ be in $M(\cA,\cW)$.
We show that the series $t_{ij}:=\sum_{k=1}^\infty r_{ik}s_{kj}$
converge absolutely in $\cA$, and that $T:=(t_{ij})_{i,j\in \N}\in M(\cA,\cW)$.
To this end, let $f,g\in\cW$ with $C_g<\infty$.
If $i\geq j$, we have $i\vee j = i \leq i\vee k$ for all $k\in \N$,
hence $f(i\vee j)\leq f(i\vee k)$ by monotonicity and thus
\begin{eqnarray}
f(i\vee j)\sum_{k=1}^\infty \| r_{ik}\|\,\|s_{kj}\| & = &
\sum_{k=1}^\infty f(i\vee j)\| r_{ik}\|\,\|s_{kj}\|
\leq \sum_{k=1}^\infty f(i\vee k)\| r_{ik}\|\,\|s_{kj}\|\notag \\
&\leq&  \|R\|_f\sum_{k=1}^\infty  \|s_{kj}\|
\leq \|R\|_f\sum_{k=1}^\infty  \underbrace{g(k\vee j)\|s_{kj}\|}_{\leq \|S\|_g}\frac{1}{g(k\vee j)}\notag \\
& \leq&  \|R\|_f\|S\|_g\sum_{k=1}^\infty \frac{1}{g(k\vee j)}\leq 
C_g\, \|R\|_f\|S\|_g < \infty,\label{esi1}
\end{eqnarray}
using that $g$ is monotonically increasing for the penultimate inequality.
If $i\leq j$, the same argument shows that
\begin{equation}\label{esi2}
f(i\vee j)\sum_{k=1}^\infty \| r_{ik}\|\,\|s_{kj}\| \leq 
C_g\, \|R\|_g\|S\|_f < \infty.
\end{equation}
In particular, in either case $\sum_{k=1}^\infty \| r_{ik}\|\,\|s_{kj}\| <\infty$, whence indeed
$\sum_{k=1}^\infty  r_{ik}s_{kj}$ converges absolutely.
Now (\ref{esi1}) and (\ref{esi2}) show that $SR:=T:=(t_{ij})_{i,j\in\N}\in M(\cA,\cW)$,
with
\begin{equation}\label{esi3}
\|SR\|_f\leq C_g\, (\|R\|_f\|S\|_g\vee \|R\|_g\|S\|_f)\,.
\end{equation}
Step 2: We show that the multiplication just defined is associative.
To this end, let $R=(r_{ij})_{i,j\in \N}$,
$S=(s_{ij})_{i,j\in \N}$ and
$T=(t_{ij})_{i,j\in \N}$ be in $M(\cA,\cW)$.
Let $R'$, $S'$ and $T'$ be the matrices
with entries $\|r_{ij}\|$, $\|s_{ij}\|$ and $\|t_{ij}\|$,
respectively.
Then $R',S',T'\in M(\R,\cW)$, as is clear from the definitions.
Hence
\begin{eqnarray*}
\sum_{(k,\ell)\in \N\times \N} \|r_{i\ell}\|\, \|s_{\ell k}\|\, \|t_{kj}\|
& =& \sum_{\ell=1}^\infty\sum_{k=1}^\infty \|r_{i\ell}\|\, \|s_{\ell k}\|\, \|t_{kj}\|
\; =\; \sum_{\ell=1}^\infty \|r_{i\ell}\|(S'T')_{\ell j}\\
& =& (R'(S'T'))_{ij}\; \in \; \R
\end{eqnarray*}
(where the first equality is a well-known elementary fact,
which can also be infered by applying
Fubini's Theorem to the counting measures on $\N^2$ and $\N$).
Thus
$\sum_{(k,\ell)\in \N\times \N} \|r_{i\ell}s_{\ell k}t_{kj}\|<\infty$,
showing that the family $(r_{ik}s_{k\ell}t_{\ell j})_{(k,\ell)\in \N\times \N}$
of elements of $\cA$ is absolutely summable.
As a consequence,
\begin{eqnarray*}
((RS)T)_{ij} &=& \sum_{k=1}^\infty (RS)_{ik}t_{kj}
=\sum_{k=1}^\infty\sum_{\ell=1}^\infty r_{i\ell}s_{\ell k}t_{kj}
=\sum_{(k,\ell)\in \N\times \N} r_{i\ell}s_{\ell k}t_{kj}\\
&=&\sum_{\ell=1}^\infty\sum_{k=1}^\infty r_{i\ell}s_{\ell k}t_{kj}
=\sum_{\ell=1}^\infty r_{i\ell}(ST)_{\ell j}
=(R(ST))_{ij}
\end{eqnarray*}
using \cite[5.3.6]{Die} for the third and fourth equalities.
Thus $(RS)T=R(ST)$.\\[3mm]
Step 3. The locally convex space $M(\cA,\cW)$ is complete.
To see this, note first that $M(\cA,\{f\})$ (with the norm $\|.\|_f$)
is a Banach space isomorphic to the space $\ell^\infty(\cA)$ of
bounded $\cA$-valued sequences, for each $f\in \cW$.
Next, after replacing $\cW$ with the set of finite sums
of elements of $\cW$ (which changes neither $M(\cA,\cW)$ as a set, nor its topology),
we may assume henceforth that $\cW+\cW\sub \cW$ and hence that $\cW$ is upward directed.
Then $M(\cA,\cW)$ is the projective limit of the
complete spaces $M(\cA,\{f\})$ ($f\in \cW$) and hence complete.\\[3mm]
Step 4. We show that the set $Q$ of quasi-invertible elements in $M(\cA,\cW)$ is open.
By
\cite[Lemma 2.6]{Glo}, we need only check that $Q$ is a $0$-neighbourhood.
To this end, choose $g\in \cW$ such that $C_g<\infty$.
Then
\[
\Big\{T\in M(\cA,\cW)\colon \|T\|_g<\frac{1}{C_g}\Big\}\;\sub \; Q\,.
\]
Indeed, pick $T$ in the left hand side.
We claim that
\begin{equation}\label{byind}
(\forall n\in \N)\qquad \|T^n\|_f\;\leq\;
(C_g\|T\|_g)^{n-1}\|T\|_f
\end{equation}
for each $f\in \cW$.
If this is true, then $\sum_{n=1}^\infty T^n$ converges in each of the Banach spaces
$(M(\cA,\{f\}),\|.\|_f)$ and hence also in the projective limit
$M(\cA,\cW)$. Now the usual argument shows that
$-\sum_{n=1}^\infty T^n$ is the quasi-inverse of $T$.\\[3mm]
To prove the claim, we proceed by induction.
If $n=1$, then $\|T\|_f=(C_g\|T\|_g)^0\|T\|_f$.
If the claim holds for $n-1$ in place of $n$, writing $T^n=T^{n-1}T$ we deduce from (\ref{esi3}) that
\begin{equation}\label{s1}
\|T^n\|_f \leq  C_g\, (\|T^{n-1}\|_f\|T\|_g\vee \|T^{n-1}\|_g\|T\|_f)\,.
\end{equation}
Now
\begin{equation}\label{s2}
C_g\|T^{n-1}\|_f\|T\|_g \leq
C_g(C_g\|T\|_g)^{n-2}\|T\|_f \|T\|_g=(C_g\|T\|_g)^{n-1}\|T\|_f
\end{equation}
by induction.
Likewise,
\begin{equation}\label{s3}
C_g\|T^{n-1}\|_g\|T\|_f\leq
C_g(C_g\|T\|_g)^{n-2}\|T\|_g
\|T\|_f=(C_g\|T\|_g)^{n-1}\|T\|_f,
\end{equation}
applying the inductive hypothesis to $g$ and $g$
in place of $f$ and $g$.
Combining (\ref{s1}),
(\ref{s2}) and (\ref{s3}),
we see that $\|T^n\|_f\leq
(C_g\|T\|_g)^{n-1}\|T\|_f$, which completes the inductive proof.\\[3mm]
Step 5. $M(\cA,\cW)$ is locally m-convex.
To see this, pick $g\in \cW$ with $C_g<\infty$.
After replacing $\cW$ with $\{f+g\colon f\in \cW\}$
(which changes neither $M(\cA,\cW)$ as a set nor its topology),
we may assume henceforth that $C_f<\infty$ for each $f\in \cW$.
We may therefore choose $g:=f$ in (\ref{esi3}) and obtain
\[
\|RS\|_f\leq C_f\|R\|_f\|S\|_f\,.
\]
Let $h:=C_f\cdot f$. Then $C_h=\frac{1}{C_f}\sum_{n=1}^\infty \frac{1}{f(n)}=1$
and $\|.\|_f$ is equivalent to the norm $\|.\|_h$,
which is submultiplicative as
$\|RS\|_h\leq C_h\|R\|_h\|S\|_h=\|R\|_h\|S\|_h$.\\[3mm]
Step 6. Continuity of quasi-inversion.
Since we assume that $C_f<\infty$ for each $f\in \cW$,
we know from Step 5 that $M(\cA,\{f\})$ is a Banach algebra,
with respect to a submultiplicative norm $\|.\|_h$ which is equivalent to $\|.\|_f$.
Now, as we assume that $\cW+\cW\sub \cW$ (see Step 3),
$M(\cA,\cW)$ is the projective limit of the Banach algebras
$M(\cA,\{f\})$ ($f\in \cW$).
Because quasi-inversion is continuous
in each of the Banach algebras,
and continuity of maps into projective limits can be checked
componentwise,
it follows that quasi-inversion is continuous on $Q\sub M(\cA,\cW)$.\,\Punkt\vspace{3mm}

\noindent
{\bf Remark.} Our results were first recorded in the unpublished thesis \cite{Lan}.\vspace{-3.5mm}
Corresponding author:\\[1mm]
Helge Gl\"{o}ckner, University of Paderborn, Institute of Mathematics,\\
Warburger Str.\ 100, 33098 Paderborn, Germany.\\[1mm]
E-Mail: {\tt  glockner\at{}math.upb.de}

\begin{thebibliography}{99}
%
%
\bibitem{Die}
Dieudonn\'{e}, J., ``Foundations of Modern Analysis,'' Academic Press, 1969.
%
%
\bibitem{Glo}
Gl\"{o}ckner, H., \emph{Algebras whose groups of units are Lie groups},
Studia Math.\ {\bf 153} (2002), no.\ 2, 147--177.
%
%
\bibitem{Lan} Langkamp, B., ``Banachalgebren und Algebren mit stetiger Inversion,''
Bachelorarbeit, Universit\"at Paderborn, 2010
(advised by  H. Gl\"ockner).
%
\bibitem{Nee} Neeb, K.-H., \emph{Towards a Lie theory of locally convex groups},
Jpn.\ J.\linebreak
Math.\ {\bf 1} (2006), No.\ 2, 291--468.
%
%
\bibitem{Mic} Michael, E., ``Locally Multiplicatively Convex Topological Algebras,'' Mem.\  Am.\ Math.\ Soc.\ 
{\bf 11}, 1952.
%
\bibitem{Phi} Phillips, N.\,C., \emph{K-theory for Fr\'{e}chet algebras},
Internat.\ J.\ Math.\ {\bf 2} (1991), 77--129.
%
%
\bibitem{Wal}
Waelbroeck, L., \emph{Les alg\`{e}bres \`{a} inverse continu},
C. R. Acad.\ Sci., Paris {\bf 238} (1954), 640--641.\vspace{1mm}
%
\end{thebibliography}
\end{document}